%%%%%%% Twisted sums and a problem of Klee
%%%%%%% by N. Tenny Peck to Victor Klee
%%%%%%% typed by Hilda K. Britt October 1991
%%%%%%%%%%%%%%%%%%%%%%%%%%%%%%%%%

\documentstyle{amsppt}
\magnification=1200
\NoBlackBoxes
\TagsOnRight

%\hsize=6.5truein
%\vsize=22.5truecm

%%\define\lrangle#1{{\langle #1 \rangle}}
%%\define\Gal{\operatorname{Gal}}
\define\a{\alpha}

\define\e{\eta}

\redefine\t{\tau}
\define\({\left(}
\define\){\right)}
\define\[{\left]}
\define\]{\right]}
\define\Span{{\operatorname{span}}}

\topmatter
\title Twisted sums and a problem of Klee
\endtitle
\author N. T. Peck \endauthor
\dedicatory To Victor Klee \enddedicatory
\address University of Illinois, Urbana, IL 61801 \endaddress
\endtopmatter

\document
\baselineskip=24pt plus 2pt

In \cite5, Klee asked whether every vector topology $\t$ on a 
real vector space $X$ is the supremum of a nearly convex topology 
$\t_1$ and a trivial dual topology $\t_2$. Recall that a vector 
topology $\t_1$ on $X$ is {\it nearly convex} if for every $x$ 
not in the $\t_1$--closure of $\{0\}$ there is $f$ in 
$(X,\t_1)^\ast$ with $f(x)\ne 0$; $\t_2$ is {\it trivial dual} if 
$(X,\t_2)^\ast=\{0\}$. We do not require that $\t_1$ or $\t_2$ be 
Hausdorff, even if $\t$ itself is Hausdorff. The topology $\t$ is 
the {\it supremum} of $\t_1$ and $\t_2$ if $\t_1$ and $\t_2$ are 
weaker than $\t$, and if for every $\t$--neighborhood $U$ of the 
origin $0$ there are a $\t_1$--neighborhood $V$ of $0$ and a 
$\t_2$--neighborhood $W$ of $0$ such that $U\supset V\cap W$.

In \cite5, Klee proved that the usual topology on $\ell_p$, 
$0<p<1$, is not the supremum of a locally convex topology and a 
trivial dual topology; this and other examples make the question 
at the beginning of this paper a natural one.  Some related 
questions on suprema of linear topologies were studied in \cite7.

Given any vector topology $\t$ on $X$, let $K(\t)=\cap \{f^{- 1}(0) 
: f\in (X,\t)^\ast\}$. It is trivial to answer Klee's question 
affirmatively in the case that $K(\t)$ is complemented. For in 
this case, $K(\t)$ must be a trivial dual space in the relative 
topology; and if $L$ is a complement to $K(\t)$ in $X$, the 
relative topology on $L$ is nearly convex. Now simply let $\t_1$ 
be the product of the trivial topology on $K(\t)$ and the 
relative topology on $L$; and let $\t_2$ be 
the product of the relative topology on $K(\t)$ and the trivial 
topology on $L$. Then $\t=\sup(\t_1,\t_2)$.

So the interesting case is when $K(\t)$ is uncomplemented. We 
study the problem when $(X,\t)$ is the {\it twisted sum} of a 
separable normed space and the real line. Recall that a real 
function $F$ on a normed space $E$ is {\it quasi--linear} if
\roster
\item"{(0)\,\,(i)}"\,\, $F(rx)=rF(x)$ for all scalars $r$ and all $x$ in 
$E$;

\item"{\phantom{(0)\,\,}(ii)}"\,\, 
$|F(x+y)-F(x)-F(y)| \leq C(\|x\| + \|y\|)$ for all $x$, $y$ 
in $E$ and some constant $C$. 
\endroster 

Now define the {\it twisted sum} of the real line and $E$ (with 
respect to $F$) as the vector space $X_F= \Bbb R \times E$ 
equipped with quasi--norm $\||(r,x)|\| = |r-F(x)| + \|x\|$. It is 
easy to verify that
$$
\||(r_1 +r_2, x_1+x_2)|\|\leq (C+1) [\||(r_1,x_1)|\| + 
\||(r_2,x_2)|\|].
$$

The space $E$ is said to be a $K$--{\it space\/} if the subspace 
$\Bbb R \times \{0\}$ is complemented in $X_F$ for every 
quasi--linear map $F$  on $E$. (This is a slight abuse of terminology; 
strictly speaking, it is the completion of $E$ that is the 
$K$--space.) So we are interested in Klee's question for the 
non--$K$ spaces. The only known non--$K$ spaces are 
$\ell_1$--like.  The {\it Ribe function\/} is defined on 
$\ell^0_1$, the space of finitely supported elements of $\ell_1$, 
by
$$
F_0(x) = \sum_i x_i \ell n|x_i| - \( \sum_i x_i \) \ell n \big|
\sum_i x_i \big|
$$
with the convention that $0\ell n 0 = 0$. Ribe \cite8 proved that 
$F_0$ is quasi--linear on $\ell^0_1$ and used $F_0$ 
to show that $\ell_1$ is not a $K$--space. Closely related 
functions were used by Kalton \cite2 and Roberts \cite9 to prove 
the same result. The reflexive space $\ell_2(\ell^n_1)$ is not a
$K$--space, and the $B$--convex spaces are $K$--spaces \cite2. 
Kalton and Roberts \cite4 showed that $c_0$ and $\ell_\infty$ are 
$K$--spaces. It is not known whether the James space is a 
$K$--space. We are studying Klee's problem for spaces $E$ and quasi-linear 
maps $F$ on $E$ such that $\Bbb R \times \{0\}$ is not 
complemented in $X_F$. By Theorem 2.5 of \cite3, there is no 
linear map $T$ on $E$ such that $|T(x)-F(x)|\leq C\|x\|$ for all 
$x$ in $E$ (i.e. $F$ does not split on $E$). The corollary to our main
theorem implies that none of the spaces above can be a counterexample
for Klee's question, since the $F$ concerned does split on an infinite-
dimensional subspace.

We now state our main result:

\proclaim{Main Theorem} Let $E$ be an $\aleph_0$--dimensional normed 
space. Assume $F$ is a quasi--linear function on $E$ for which 
there are a linearly independent sequence $(x_i)$ in $E$ and a 
linear map $T$ on $\Span(x_i)$ such that
\roster
\item"{(1)}" $|T(x)-F(x)| \leq C\|x\|$ for all $x$ in 
$\Span(x_i)$ and some constant $C$.
\endroster

\noindent
Then there are a trivial dual topology $\t_2$ on $\Bbb R\times 
E$, weaker than the quasi--norm topology, and a 
$\t_2$--neighborhood $U$ of $0$ such that if $(r,x) \in U$ and 
$\|x\| \leq 1$, then $\||(r,x)|\|<C$ for some constant $C$.
\endproclaim

Before we prove the theorem, we set the framework for the 
construction with some auxiliary results. We begin with: 

\definition{Definition} Suppose $(G_i)$ is a finite or infinite 
sequence of subsets of $E$, and $(n_i)$ is a sequence of positive 
integers (of the same length as $(G_i)$). The $(n_i)$--sum of 
$(G_i)$ is the set of all finite sums
$$
z=r_1z_1 + r_2z_2 + r_3z_3\ldots
$$
where $|r_i| \leq 1$ for all $i$ and $z_1,\ldots,z_{n_1}$ are in 
$G_1$, $z_{n_1+1},\ldots,z_{n_1+n_2}$ are in $G_2$,\newline 
$z_{n_1+n_2+1},\ldots,z_{n_1+n_2+n_3}$ are in $G_3$, etc. Note 
that if $|r|\leq 1$, $rz$ is also in the $(n_i)$--sum.
\enddefinition

\proclaim{Lemma 1} Let $X$ be a vector space and let $(U_n)$ be 
a neighborhood base at $0$ for a pseudo-metrizable vector topology on 
$X$, chosen so that $U_{n+1} + U_{n+1} \subset U_n$ for all $n$ 
and $[-1,1]U_n\subset U_n$ for all $n$. Let $(F_n)$ be a sequence 
of subsets of $X$, chosen so that $[-1,1]F_n\subset F_n$ and 
$F_{n+1} + F_{n+1} \subset F_n$, for all $n$. Then the sequence $(U_n + F_n)$ 
is a neighborhood base at $0$ for a pseudo-metrizable vector topology on 
$X$ which is weaker than the original topology.
\endproclaim

\demo{Proof} Immediate. \qed
\enddemo

In the next lemma, we specify $F_n$ more closely.

\proclaim{Lemma 2} Let $X$ and $(U_n)$ be as in Lemma 1. Let 
$(G_n)$ be a sequence of subsets of $X$. Define subsets $F_n$ of 
$X$ as follows:  for each $n$ in $N$, $F_n$ is the 
$(2^{i-n})$--sum of the $G_i$'s for $i\geq n$. 
Then $(F_n)$ satisfies the hypotheses of Lemma 1.
\endproclaim

\demo{Proof} $[-1,1]F_n \subset F_n$ as remarked already. For a 
typical sum in $F_{n+1} + F_{n+1}$, at most $2\cdot 2^{i-(n+1)} = 
2^{i-n}$ of the $z_i$'s are in $G_i$ for $i\geq n+1$, so $F_{n+1} 
+ F_{n+1} \subset F_n$. \qed
\enddemo

\remark{Remark} Note that there is an apriori bound on the number 
of elements of $G_i$ appearing in a sum in $F_n$, for any $n$: 
the bound is $2^{i-1}$; we use the looser bound $2^i$.

In our construction, $(U_n)$ is a neighborhood base at $0$ for 
the twisted sum topology. The $G_i$'s of Lemma 2 will be chosen 
so that $(U_n+F_n)$ is a neighborhood base at $0$ for a trivial 
dual topology $\t_2$; they will also have to be chosen so that 
$\t$ is the supremum of $\t_1$ and $\t_2$. The next lemma 
identifies the topology $\t_1$:
\endremark

\proclaim{Lemma 3} Let $F$ be a quasi--linear map on a normed 
space $E$ and let $X_F=\Bbb R \times E$ with the quasi--norm 
$\||(r,x)|\| = |r-F(x)| + \|x\|$. Assume $\Bbb R\times \{0\}$ 
is not complemented in $X_F$. Then the strongest nearly convex 
topology on $\Bbb R\times E$ which is weaker than the 
quasi--norm topology has a neighborhood base at $0$ of sets of 
the form $\{ (r,x) : \|x\| < \e\}$.
\endproclaim

\demo{Proof} Sets of the above type are a neighborhood base at 
$0$ for a nearly convex topology weaker than $\t$, the quasi-norm 
topology. The closure of $\{0\}$ for this weaker topology is 
$\Bbb R \times \{0\}$; and if $(r,x)$ is in $X_F$ and $x\ne 0$, 
there is $f$ in $E^\ast$ with $f(x)\ne 0$. Then $f(\pi(r,x)) \ne 
0$, where $\pi$ is the quotient map of $X_F$ onto $E$.

Now suppose $\nu$ is a nearly convex topology on $\Bbb R \times 
E$, weaker than the quasi-norm topology. Since $\Bbb R \times 
\{0\}=K(\tau)$ is not complemented, the $\nu$--closure of $\{0\}$ 
must contain $\Bbb R \times \{0\}$. Let $U$ be $\nu$--open 
containing $0$. Choose $V$ $\nu$--open containing $0$, with $V+V\subset 
U$. Choose $\epsilon>0$ so that if $\|| (r,x) |\| <\epsilon$, then 
$(r,x)\in V$.

Now, $\||(F(x),x)|\| = \| x\|$, so if $\|x\| <\epsilon$, then 
$(F(x),x)\in V$. Also, $(r-F(x),0)$ is in $V$ since it is in the 
$\nu$--closure of $0$, and so $(r,x)$ is in $U$. \qed
\enddemo

\definition{Notation} Let $Z$ be a Banach space with a basis 
$(v_i)$. Let $(v^\ast_i)$ be the coordinate functionals on $Z$. 
For $n$ a positive integer and $x$ in $Z$, set $x\mid_{[1,n]} = 
\sum\limits^n_{i=1} v^\ast_i(x)v_i$, and  $x\mid_{(n,\infty)} = 
\sum\limits^\infty_{i=n+1} v^\ast_i(x)v_i$. Say that $x$ {\it is 
to the right of} $n$ if $v^\ast_i(x)=0$ for $i\leq n$.
\enddefinition

We need two more preliminary results before proving our main 
theorem:

\proclaim{Lemma 4} Let $Z$ be a Banach space with a monotone 
basis $(v_i)$, let $K$ be a compact subset of $Z$, and let 
$\epsilon >0$. Then there is $n$ so that if $y$ is to the right 
of $n$ and $x\in K$, then $\| x\| < \|x+y\| + \epsilon$.
\endproclaim

\demo{Proof} Choose $n$ so that $\|x \mid_{(n,\infty)} \| 
<\epsilon$ for every $x$ in $K$. Now if $y$ is to the right of 
$n$, then $\| x \mid_{[1,n]}\|= \|(x+y)\mid_{[1,n]}\| \leq 
\|(x+y)\|$, since the basis is monotone, so $\|x\| < \|x+y\| + 
\epsilon$. \qed 
\enddemo 

\proclaim{Lemma 5} Let $y_i$, $1\leq i \leq k$, be linearly 
independent elements of a normed space $E$. Define $y_{k+1} =-
\sum\limits^k_{i=1} y_i$, and let $\eta >0$. Set $z_i = my_i$, 
$1\leq i\leq k+1$, for some $m>0$. Then  we can choose $m$ so 
large that the following condition is satisfied: if at most $k$ of 
the $r_i$'s are non--zero, and if $\big\| \sum\limits^{k+1}_{i=1} 
r_iz_i\big\| < 3$, then $\sum\limits^{k+1}_{i=1} |r_i| < \eta$.
\endproclaim

\demo{Proof} Choose $M>0$ so that $\sum\limits^k_{i=1} |\a_i| \leq 
M \big\| \sum\limits^k_{i=1}\a_iy_i\big\|$ for all $k$--tuples $(\a_i)$. 
Now suppose  $\big\| \sum\limits^{k+1}_{i=1} r_iz_i\big\| < 3$, with at 
most $k$ $r_i$'s non--zero. If $r_{k+1}=0$, then 
$$
\sum^k_{i=1} |r_i| \leq \frac{3M}m < \eta
$$
for $m>3M/\eta$.

If $r_{k+1}\ne 0$, then some other $r_i$ is $0$, $r_1$, say; now,
$$
\big\| \sum^{k+1}_{i=1} r_iz_i\big\| = \big\| -m r_{k+1} y_1 + \sum^k_{i=2} 
m(r_i-r_{k+1})y_i\big\| <3.
$$
Therefore
$$
\aligned
|r_{k+1}| 
&+ \sum^{k+1}_{i=2} |r_i-r_{k+1}| < \frac{3M}m, \quad \text{so} \\
   \sum^{k+1}_{i=2} |r_i| 
&\leq |r_{k+1}| + \sum^{k+1}_{i=2} |r_i-r_{k+1}| + k|r_{k+1}| \\
&< \frac{3(k+1)M}m <\eta
\endaligned
$$
for $m>3(k+1)M/\eta$. \qed
\enddemo

\demo{Proof of Main Theorem} We will construct inductively the 
sets $G_n$ used in Lemma 2. That lemma will give us the sets 
$F_n$ and then Lemma 1 will provide the topology. 

We may assume that for $x_i$ and $T$ in the theorem, $T(x_i) = 0$ 
for each $i$. This is possible since for each $i$, there is a 
scalar $\a_i$ so that $T(x_{2i-1} + \a_ix_{2i}) =0$; now the 
sequence $x'_i=x_{2i-1} + \a_i x_{2i}$ also satisfies condition 
1) of the Theorem.

We can regard $E$ as a subspace of the Banach space $Z=C[0,1]$, 
which has a monotone basis.  Any positive scalar 
multiple of the quasi--norm yields the same topology as the 
quasi--norm, so we can and do assume that the constant $C$ 
in $0$(ii) above is 
$1$. This can be done by multiplying $F$ by a suitable positive 
constant.

Finally, we use $\| \quad \|$ to refer to the norm on $Z$ and on 
$E$. We only calculate norms of  elements of $E$, but we do use 
the monotonicity of $(v_i)$ in $Z$.

Now we begin the construction of $(G_i)$.

Choose $0<c_n\leq 2^{-(n+3)}$ (and thus $\sum\limits^\infty_{n=1} 
c_n < \dfrac14$). Let $(d_j)$ be any sequence whose linear span is 
$E$, and let $(e_i)$ be an indexing of $(d_j)$ such that each $d_j$ occurs 
infinitely often in $(e_i)$. We can assume  that $\|d_j\| \leq 1$ 
and that $|F(d_j)|\leq 1$ for each $j$, by multiplying $d_j$ by 
a positive constant.

Assume that finite sets $G_0$, $G_1\cdots G_{n-1}$ have been 
constructed, with $G_0=\{0\}$, satisfying the following conditions:
\roster
\item"{(2)}" for each $1\leq i\leq n-1$, $G_i$ is a finite set 
$(w_{i,j} : 1\leq j \leq 2^{i+1})$, with
$$
\align
w_{i,j}     &= e_i + m_i x_{\ell(i,j)}, \quad j \leq 2^i, \\
w_{i,2^i+1} &= e_i - \sum^{2^i}_{j=1} m_i x_{\ell(i,j)}; 
\endalign
$$
here, $(x_i)$ is the sequence in the statement of the theorem.

\item"{(3)}" Set $z_{i,j} = w_{i,j} - e_i$. Then if 
$\big\|\sum\limits^{2^i+1}_{j=1} r_j z_{i,j} \big\| < 3$ with at most 
$2^i$ $r_j$'s non--zero, $\sum\limits^{2^i+1}_{j=1} |r_j|<c_i$.
\endroster

To define $G_n$, let $K_n'$ be the $(2^i)$--sum of $G_i$ for 
$i\leq n-1$, (so $K_1'=\{0\}$) and let $K_n = K_n' + [-
2^n,2^n]e_n$. Then $K_n$ is a compact subset of $E\subset Z$. By 
Lemma 4, there is an integer $s_n$ such that if $y$ is to the 
right of $s_n$ then $\|x\| < \| x+y\| + c_n$ for all $x$ in 
$K_n$. By the linear independence of the sequence $(x_i)$, we can 
choose $x_{\ell(n,1)},\ldots, x_{\ell(n,2^n)}$, all to the right 
of $s_n$, with $\ell(n,j)<\ell(n,j')$ if $j<j'$. For ease of 
notation, put $x_{n,i} = x_{\ell(n,i)}$, $1\leq i\leq 2^n$, and 
put $x_{n,2^n+1} = - \sum\limits^{2^n}_{i=1} x_{n,i}$.

By Lemma 5, we can choose $m_n$ so large that if
$$
\big\| \sum\limits^{2^n+1}_{j=1} r_{n,j} m_n x_{n,j}\big\| <3,
$$
with at most $2^n$ of the $r_{n,j}$ non--zero, then
$$
\sum\limits^{2^n+1}_{j=1} |r_{n,j}|<c_n.
\tag4
$$
Finally, for $1\leq i \leq 2^n+1$, put
$$
w_{n,i} = e_n + m_nx_{n,i}
$$
and let
$$
G_n = (w_{n,i}), \quad 1\leq i \leq 2^n +1.
$$
Note that since $\sum\limits^{2^n+1}_{i=1} x_{n,i} =0$, $e_n\in 
co G_n$. (We denote the convex hull of $A$ by $co A$.) This 
finishes the construction of $(G_i)$.

Now let $(F_n)$ be the subsets of $E$ used in Lemma 1: $F_n$ is 
the $(2^{i-n})$ sum of $(G_i)$ for $i\geq n$. Let $(U_n)$ be a 
neighborhood base at $0$ for the quasi--norm topology on $\Bbb 
R\times E$, with $U_{n+1} + U_{n+1} \subset U_n$ and $[-
1,1]U_n\subset U_n$, for all $n$; also assume that $\||w|\| <1$ if 
$w\in U_1$. Let $\t_2$ be the topology yielded by Lemma 1.

We claim that $\t_2$ is trivial dual. To see this, note that for 
$m\geq n$, $e_m\in co(w_{m,i}) \subset co F_n \subset co(F_n + 
U_n)$; since each $d_j$ occurs infinitely often in the sequence 
$(e_m)$, $K(\t_2)$ contains every $d_j$ and therefore contains 
$\{0\}\times E$. Also, $(1,0)\in co U_n \subset co(U_n+F_n)$ for every $n$, 
so $K(\t_2)$ contains $\Bbb R\times \{0\}$. This proves the claim.

Now suppose that
$$
x=\sum^n_{i=1} \sum^{2^i+1}_{j=1} r_{i,j}(e_i + m_ix_{i,j})
$$
is in $F_1$ and that $\|x\|<1$. We will first prove that 
$|F(x)|<9$. 

Toward that end: since the $x_{n,j}$ are to the right of $s_n$, the 
construction of $G_n$ implies that
$$
\big\| \sum^{n-1}_{i=1} \sum^{2^i+1}_{j=1} r_{i,j}(e_i + m_i x_{i,j}) 
+ \sum^{2^n+1}_{j=1} r_{n,j} e_n\big\| <1+ c_n
\tag5
$$
from which, since $\|x\|<1$,
$$
\big\| \sum^{2^n+1}_{j=1} r_{n,j}m_nx_{n,j}\big\| <2 + c_n<3.
\tag6
$$
Now from (4) and (6), we have
$$
\sum^{2^n+1}_{j=1}|r_{n,j}|<c_n;
\tag7
$$
combining this with (5), we have
$$
\big\|\sum^{n-1}_{i=1} \sum^{2^i+1}_{j=1} r_{i,j}(e_i + m_i 
x_{i,j})\big\| <1+2c_n.
\tag8
$$
For the induction step, assume that for some $\ell$,
$$
\big\|\sum^\ell_{i=1} \sum^{2^i+1}_{j=1} r_{i,j}(e_i + m_i 
x_{i,j})\big\| < 1 + 2c_n \cdots+ 2c_{\ell+1}.
\tag9
$$
Since the $x_{\ell,j}$ are to the right of $s_\ell$, the 
construction of $G_\ell$ implies that
$$
\big\| \sum^{\ell-1}_{i=1} \sum^{2^i+1}_{j=1} r_{i,j}(e_i + m_i 
x_{i,j}) + \sum^{2^\ell+1}_{j=1} r_{\ell,j}e_\ell\big\|<1+2c_n\cdots+ 
2c_{\ell+1} + c_\ell,
\tag10
$$
from which
$$
\big\|\sum^{2^\ell+1}_{j=1}r_{\ell,j}m_\ell x_{\ell,j}\big\| <2+4c_n 
\cdots + 4c_{\ell+1} + c_\ell <3.
\tag11
$$
Now from (4) and (11), we have
$$
\sum^{2^\ell+1}_{j=1}|r_{\ell,j}|<c_\ell;
\tag12
$$
combining this with (10), we obtain 
$$
\big\|\sum^{\ell-1}_{i=1} 
\sum^{2^i+1}_{j=1} r_{i,j}(e_i + m_i x_{i,j}) \big\| < 1 + 2c_n 
\cdots + 2c_\ell, 
\tag13 
$$ 
recalling that $\|e_i\|\leq 1$. This 
finishes the induction step. 

The above argument has yielded that
$$
\big\|\sum^{2^i+1}_{j=1} r_{i,j}e_i\big\|<c_i
\tag14
$$
for each $i$; from this and $\|x\|<1$, we have
$$
\big\|\sum^n_{i=1} \sum^{2^i+1}_{j=1} r_{i,j}m_ix_{i,j} \big\| <1 + 
\sum^\infty_{n=1} c_n<2.
\tag15
$$
>From (15) and (1), recalling $T(x_{i,j})=0$,
$$
\big|F\( \sum^n_{i=1} \sum^{2^i+1}_{j=1} r_{i,j}m_ix_{i,j} 
\)\big|<2.
$$
To estimate
$$
F\( \sum^n_{i=1} \sum^{2^i+1}_{j=1} r_{i,j}e_i\),
$$
recall that $|F(e_i)|\leq 1$ for each $i$, so 
$$
\big|F\(\sum^{2^i+1}_{j=1} r_{i,j}e_i\)\big|<2^{-i}.
$$
Therefore
$$
\align
\big|F\(\sum^n_{i=1}\(\sum^{2^i+1}_{j=1} r_{i,j}e_i\)\)\big| &\leq 
\sum^n_{i=1} 2^{-i} + \sum^n_{i=1} 
i\big\|\sum^{2^i+1}_{j=1}r_{i,j}e_i\big\| \\
&<1 + \sum^n_{i=1} i\cdot 2^{-i}<4
\endalign
$$
(using $\big|F(\sum u_i)\big|\leq \sum |F(u_i)| + \sum i \| u_i\|)$. 
Finally, 
$$
\aligned
|F(x)| \leq& 2 + 4 + \big\| \sum^n_{i=1} \sum^{2^i+1}_{j=1} 
r_{i,j}e_i\big\| \\
&+\big\|\sum^n_{i=1} \sum^{2^i+1}_{j=1} r_{i,j}m_ix_{i,j}\big\|\\
<& 2+4+1+2=9.
\endaligned
$$

To complete the proof of the theorem, suppose $(r,x)\in U_1 + 
F_1$ and $\|x\| \leq 1$. Write $(r,x)=(r,y) + (0,z)$, with 
$(r,y)\in U_1$ and $z\in F_1.$ Then $|r-F(y)| + \|y\|\leq 1$; 
from this and $\|x\|\leq 1$ follows $\|z\|\leq 2.$ Now since 
$z\in F_1$, the preceding paragraph implies $|F(z)|<18$. At last,
$$
\aligned
|r-F(x)| &\leq |r-F(y)| + |F(y)-F(x)| \\
&\leq 1+|F(y)-F(x)| \\
&\leq 1+ |F(z)| + \|z\| + \|x\| \\
&<22,
\endaligned
$$
so $\||(r,x)|\|<23$. The proof is complete. \qed
\enddemo

\proclaim{Corollary} Let $E$ be a separable normed space
and let $E_0$ be an $\aleph_0$-dimensional subspace of
$E$ which is dense in $E$.
Assume that there is a quasi-linear map $F$ on $E_0$ 
which splits on an infinite-dimensional subspace of
$E_0$.
Then the twisted sum topology on $\Bbb R \otimes_F E$
is the supremum of a trivial dual topology and a
nearly exotic topology. \endproclaim

\demo{Proof} Let $q$ denote the quotient map of
$\Bbb R\otimes_F \tilde E$ onto $\tilde E$, where $\tilde E$ is the 
completion of $E$. (For $x\in E_0, q(r,x)=x$.)
The subspace $E_0$ satisfies the
hypotheses of the main theorem.  Therefore there are
a trivial dual topology $\tau_2$ on $\Bbb R\times E_0$,
weaker than the twisted sum topology; a $\tau_2$-
neighborhood $V$ of $0$; and a constant $C$ so that if
$x \in E_0, (r,x) \in V$, and $\Vert x\Vert < 1$,
then $\Vert\vert (r,x)\vert\Vert < C.$

We can assume that $V$ contains a $\tau_2$-neighborhood
$U$ of $0$ of the form $B_{\alpha}+F_n$, where
$F_n\subset E_0$ is as constructed as in the proof of
the main theorem, and for any $\beta > 0$,
$$B_{\beta}=\{(r,x)\in\Bbb R\times E_0:
\Vert\vert(r,x)\vert\Vert<\beta\}.$$
Sets of the form $\overline{B_{\beta}}+q^{-1}(F_m),$ where
$$\overline{B_{\beta}}=\{w\in\Bbb R\otimes_F E:
\Vert\vert w\vert\Vert<\beta\},$$
obviously form a neighborhood base at the origin for
a vector topology $\overline{\tau}_2$ on $\Bbb R\otimes_F E,$
weaker than the twisted sum topology.  The topology
$\overline{\tau}_2$ is trivial dual since its
restriction to the dense subspace $\Bbb R\times E_0$
is trivial dual.

Now choose $0<\gamma<1/2$ so that $\overline{B_{\gamma}}+\overline{B_{\gamma}}
\subset\overline{B_{\alpha}}$, and assume that
$w\in\overline{B_{\gamma}}+q^{-1}(F_n)$ and 
$\Vert q(w)\Vert<1/2$. Choose $w_0\in
\Bbb R\times E_0$ so that $\Vert\vert w-w_0\vert\Vert
<\gamma$. Then $\Vert q(w)-q(w_0)\Vert<\gamma$,
so $\Vert q(w_0)\Vert<\gamma+1/2<1$. Clearly,
$w_0\in\overline{B_{\alpha}}+q^{-1}(F_n)$, and
so from our assumption, $\Vert\vert w_0\vert\Vert<C$.
Now, $\Vert\vert w\vert\Vert<(\alpha/\gamma)(C+1)$, and the
proof is complete. \qed
\enddemo

The theorem and corollary apply to several spaces which are either not 
$K$--spaces or for which it is not known whether they are 
$K$--spaces:

\proclaim{Theorem} For the following pairs of normed spaces $E$ 
and quasi--linear maps $F$ on $E$, the twisted sum topology on 
$X_F=\Bbb R\times E$ is the supremum of a nearly convex topology 
and a trivial dual topology:
\roster
\item"{(a)}" $E$ is any infinite--dimensional subspace of 
$\ell^0_1$ (whether or not it is a $K$--space), $F$ is the Ribe 
function $F_0$;

\item"{(b)}" $E$ is the linear span of the usual unit vector 
basis for the James space, under the James norm; $F$ is any 
quasi--linear function on $E$;

\item"{(c)}" $E$ is the span of the usual unit vector basis in 
$\ell_p(\ell^n_1)$, for $1<p<\infty$ (this is a reflexive 
non--$K$ space); $F$ will be described below.
\endroster
\endproclaim

\demo{Proof} For (a), let $H=\{x\in E: \sum\limits_i x_i = 0\}$. 
Note that if $x$, $y\in H$ and $x$ and $y$ have disjoint 
supports, $F_0(x+y)=F_0(x)+F_0(y)$. Since $H$ has codimension at 
most $1$ in $E$ and $E$ is infinite dimensional, there is a 
sequence of non--zero elements $(x_i)$ in $H$ satisfying $\sup$ 
(support $x_i$) $<\inf$ (support $x_{i+1}$) for all $i$. 

As remarked above, $F_0$ is linear on $\Span(x_i)$, so if we 
define $T(x_i)=F(x_i)$, the linear function $T$ certainly 
satisfies hypothesis (1) of the theorem. Therefore the theorem 
applies to $E$.  

For (b), it is known that the even unit vectors 
$e_{2n}$ span a pre--Hilbert subspace of the James space (see 
\cite1). The $B$--convexity of $\Span(e_{2n})$ and Theorems 2.6 
of \cite2 and 2.5 of \cite3 imply that there is a linear map $T$ 
on $\Span(e_{2n})$ such that $|T(x)-F(x)| \leq C\|x\|$ for all 
$x$  in $\Span(e_{2n})$. Therefore the theorem applies.

(c)  For each $n$ let $(e_{i,n})$ be the usual unit vector basis 
of $\ell^n_1$, and let $E$ be the span of the $e_{i,n}$ in 
$\ell_p(\ell^n_1)$. Let $(c_n)$ be any sequence in $\ell_q$, 
$\frac1p  + \frac1q = 1$. Let $F_0$ be the Ribe function and 
define $F$ on $E$ by
$$
F((x_n)) = \sum_n c_nF_0(x_n).
$$

We claim that $F$ is quasi--linear. For this, if $(x_n)$ and 
$(y_n)$ are in $E$, the sequences $(\|x_n\|_1)$ and $(\|y_n\|_1)$ 
are $\ell_p$ sequences, and for each $n$,
$$
\align
|c_n F_0(x_n + y_n) &- c_n F_0(x_n) - c_n F_0(y_n)| \\
\leq& c_n(\|x_n\|_1 + \|y_n \|_1).
\endalign
$$
>From H\"older's inequality,
$$
\align
|F((x_n + y_n)) &- F((x_n)) - F((y_n))| \\
\leq& \|(c_n)\|_q (\|(x_n)\|_p + \|(x_n)\|_p).
\endalign
$$

Theorem 4.7 of \cite2 gives that $E$ is not a $K$--space. The $F$ 
just defined proves this directly, for suppose there is a linear 
$T$ on $E$ with $|T(x)-F(x)|\leq C\|x\|$ for all $x$ in $E$. Then 
since $F(e_{i,n})=0$ for all $i$, $n$, $|T(e_{i,n})|\leq C$ for 
all $i$, $n$.  But $F\(\dfrac1n \sum\limits^n_{i=1} e_{i,n}\) = - 
c_n \log n$, a contradiction if we choose $c_n$ so that $(c_n\log 
n)$ is unbounded.

Finally, our theorem applies in this situation. To show this, for 
each $n$ pick a unit vector $x_n$ in $\ell^n_1$. The sequence 
$(x_n)$ is equivalent to the usual basis of $\ell_p$, which is 
$B$--convex; the results already mentioned imply that there is a 
linear $T$ on $\Span(x_n)$ such that $|T(x)-F(x)|\leq C\|x\|$ for 
all $x$ in $\Span(x_n)$. This finishes the proof. \qed
\enddemo

Note that, because of the separability, the corollary applies to
the completions of the twisted sums in (a)--(c) above.

\enddocument

\enddocument